\documentclass{article}
\pdfoutput=1

\usepackage{arxiv}

\usepackage[utf8]{inputenc} 
\usepackage[T1]{fontenc}    
\usepackage{hyperref}       
\usepackage{url}            
\usepackage{booktabs}       
\usepackage{amsfonts}       
\usepackage{amssymb}
\usepackage{nicefrac}       
\usepackage{microtype}      
\usepackage{lipsum}			
\usepackage{graphicx}
\usepackage{natbib}
\usepackage{doi}
\usepackage{booktabs}
\usepackage{amsmath}
\usepackage{svg}

\title{Madhava's Pi-series in the modern context}

\author{ \hspace{1mm}Harikrishna VJ \\
	Mysuru\\
	India\\
	\texttt{harikrishna.vj@gmail.com} \\
	 \and
	 \hspace{1mm}\textbf{Vittal Rao} \\
	 Bengaluru \\
	 India \\
	\texttt{vittalrao.14@gmail.com} \\
	\and
	\textbf{Srikrishna Bhat} \\
	Bengaluru \\
	India\\
	\texttt{skbhat@gmail.com} \\
}
\date{}

\renewcommand{\shorttitle}{\textit{arXiv} Template}

\hypersetup{
pdftitle={A template for the arxiv style},
pdfsubject={q-bio.NC, q-bio.QM},
pdfauthor={David S.~Hippocampus, Elias D.~Striatum},
pdfkeywords={First keyword, Second keyword, More},
}

\begin{document}
\maketitle

\begin{abstract}
The $\pi$-series is attributed to Madhava, Gregory and Leibniz based on the chronology of its discovery. While this acknowledges
the fact that Madhava of \emph{sangamagrama} an Indian mathematician who lived in Kerala in the 14th century along the banks of river \emph{nila}, was the first to discover the series, the modern day proof for the $\pi$-series is based on the outline of the proof given by the famous German mathematician  Gottfried Wilhelm  Leibniz in the $17^{th}$ century, using concepts and 
rigorous definitions of Calculus many of which were developed after his lifetime. We extend the same benefit to Madhava's ideas 
and present Madhava's \emph{upapatti} i.e proof for the $\pi$-series as given in Jyesthadeva's \emph{yuktibhasha} using modern notation and definitions. We elaborate Madhava's proof at three places (a) to prove that the arc-bit sum is same as the sum of arc-bit approximations (b) to prove convergence and (c) to concisely prove by iteration the power sum approximation. Madhava's method with appropriate changes can be applied to derive the
sine and cosine series as well. Madhava also provides the correction terms to accelarate the convergence of the $\pi$-series. We compare the two proofs i.e
Madhava's and the popular proof for $\pi$ which is based on Leibniz's ideas and observe that they are different since Madhava's is based on the length of the arc where as Leibniz's relies on the area under the arc and in terms of the number of concepts from advanced calculus used in the the two proofs Madhava's proof stands out for its simplicity, thus making a case for it to be appreciated popularized and taught as the original proof of the $\pi$-series.    
\end{abstract}
\section{Introduction}
\label{sec:intro}
Beginning with the vedic \emph{sulva sutras} India has a great tradition of mathematics enriched by mathematicians like Pingala, Aryabhata, Brahmagupta, Mahavira, Bhaskara, Narayana Pandita, Madhava and Ramanujan who worked on a wide variety of problems spanning Number theory, Geometry, Algebra and Calculus with applications in Astronomy, Geography, Architecture, Trade \& Commerce and Linguistics \cite{srinivasiengar1967history}.
Europe has a rich tradition of Mathematicians dating back to the
times of Pythagoras, Eucild, Archimedes, Fibonacci before
Rennaissance and lead by Descartes, Leibniz, Newton, Euler and Gauss,
during and after Rennaissance who shaped the growth of all the branches of Modern Mathematics \cite{boyer1968history}.

The subject of this paper is the $\pi$-series which is referred to as the Madhava-Gregory-Leibniz series based on the chronology of its discovery. While this acknowledges Madhava of \emph{sangamagrama}\footnote{The word \emph{sangama} indicates that Madhava lived at a confluence of rivers possibly \emph{nila} and one of its tributaries in a place like Kudalur in Kerala or that he may have come from a place where there was a \emph{sangama} i.e a confluence of rivers possibly from Tulu speaking parts of southern Karnataka. His works and the works of his lineage of disciples were developed in places close to the banks of river \emph{nila} \cite{Divakaran2007FirstTextbook}\cite{Krishnachandran2024Madhava}.}, the Indian mathematician who lived in Kerala in the 14th century along the banks of river \emph{nila}, as the first to discover the series, the modern day proof of the $\pi$-series is based on the outline of the proof given by the famous German mathematician Gottfried Wilhelm Leibniz who along with Isaac Newton in the $17^{th}$ century contributed immensely to the development of modern Calculus. While Leibniz's ideas found the language of modern calculus with rigorous definitions many of which came after him to become the popular modern proof for the $\pi$-series the same benifit has not been extended to Madhava's proof. To address this we present Madhava's proof given in
Jyesthadeva's \emph{yuktibhasha} using modern notation and definitions of convergence. We elaborate Madhava's proof at three places to which we have assigned names from the original text as: 
\begin{enumerate}
	\item \emph{capikaranam, Making the Arc}: To prove that the arc bit sum is same as the sum of arc-bit approximations.
	\item \emph{phala parampara, Series of Results}: To prove
	convergence of the series.
	\item \emph{samaghata samkalita, Equal Power Sum}: To concisely prove by induction the power sum approximation.
\end{enumerate}
We compare the two proofs i.e
Madhava's and the popular proof based on Leibniz's ideas and observe that they are different since Madhava's is based on the length of the arc where as Leibniz's relies on the area under the arc. Moreover in terms of the number of concepts from advanced calculus used in the two proofs Madhava's proof stands out for its simplicity thus making a case for it to be appreciated popularized and taught as the original proof of the $\pi$-series.

The proof naturally extends to give the arctan series. Madhava used the same method involving summing large number of bits, iterative refinement and approximations for power sums to discover series expansions for sine and cosine (versine) series as well. We compare the steps involved in deriving the $\pi$-series with that of the sine series. Another attractive feature of Madhava's work which has received more attention than his proof is his correction strategy to accelarate the convergence of the series. This is a distinguishing feature of his work compared to the discoveries that were made later. 
\begin{figure}
	\centering
	\includegraphics[width=130mm,height=45mm]{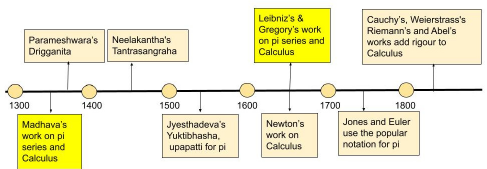}
	\caption{Timeline of development of concepts to prove the $\pi$-series }
	\label{fig:fig_1}
\end{figure}
For a timeline of the development of ideas by Madhava's school and that of European mathematicians later on please refer Figure \ref{fig:fig_1}. For more on the contributions of the European mathematicians mentioned in the timeline please refer \cite{boyer1959calculus}.

Although the work that is presented here is based on the
translation of Jyesthadeva's \emph{yuktibhasha} written in the $16^{th}$ century we refer to it as Madhava's proof since \emph{yuktibhasha} refers to \emph{tantrasangraha} written in Sanskrit by Neelakantha which attributes the results to Madhava. Also in \emph{yuktibhasha} the arctan series is explicitly attributed to Madhava \cite{Divakaran2007FirstTextbook}. Both Jyesthadeva and Neelakantha were disciples of Damodara who was a disciple of Parameshwara a direct disciple of Madhava\cite{Divakaran2007FirstTextbook}.
 
The first person to present Madhava's works to the modern world
was the English civil servant C.S Whish \cite{whish1834hindu} who worked in the Malabar region and wrote about the ``Hindu Quadrature of the Circle". Later on many mathematicians like Jushkevich \cite{Jushkevich1961Russian}, Rajagopal \cite{Rajagopal1986-RAJOMK} rediscovered these works and made them popular. The contribution of K.V Sharma \cite{sarma1972history} merits a special mention as he authored a translation of \emph{yuktibhasha} making many subsequent works like those of Ranjan Roy \cite{Roy1990Discovery}, C.K Raju \cite{Raju2002Interactions}, Plofker \cite{plofker2009mathematics}, Divakaran \cite{Divakaran2007FirstTextbook}, Joseph \cite{joseph2009passage}, M.D Srinivasan, K.Ramsubramanian and M.S Sriram possible \cite{Jyesthadeva2008}, \cite{Ramasubramanian2010}.  

There are works before ours which describe Madhava's work
by using modern notation. Among these the works of
Ranjan Roy \cite{Roy1990Discovery} and K.V Sarma, M.D Srinivasan and K. Ramsubramanian et al \cite{Jyesthadeva2008}, \cite{Ramasubramanian2010} are exemplary. But even their works do not provide a 
complete proof for the series following Madhava's thread as in Ranjan Roy's work, about Madhava's sum of arc-bit approximations which links the geometry of the formulation to algebra it is said that ``the relation was understood in an intuitive sense only" \cite[Section~4]{Roy1990Discovery}. The eventual convergence of the sum is not dealt with in detail and so is the approximation of the sum of powers. ``The Explanatory Notes" provided by  K.Ramasubramanian, M.D Srinivasan and M.S. Sriram  in ``Ganita-Yuktibhasha (Rationales in Mathematical Astronomy) of Jyesthadeva" for K.V Sarma's translation of \emph{yuktibhasha} from Malayalam \cite{Jyesthadeva2008}, forms the basis of our work. The work by K.Ramasubramanian et al \cite{Jyesthadeva2008} is also published as a paper titled `Development of Calculus in India' \cite{Ramasubramanian2010}. We retain the same diagrams and notation wherever possible to ensure continuity with the orignal translation. 

In the next section i.e Section \ref{sec:madhava} we present Madhava's \emph{upapatti} i.e proof for the $\pi$-series and the supporting Lemmas. In Section \ref{sec:madhava_other} we briefly discuss the application of Madhava's method to derive the sine series and compare the steps. In Section \ref{sec:leibniz} we present the popular modern proof of the $\pi$-series, describe Leibniz's ideas of `Transmutation' and `Integration by Parts' which eventually led to the first step of the modern proof based on the Fundamental Theorem of Calculus. In Section \ref{sec:comp} we compare the two proofs in terms of concepts from advanced calculus that go into making of the proofs. We discuss  the role played by the trigonometric inequality i.e $tan(x)>=x$, $0< x < \frac{\pi}{2}$ in the two proofs, the power sum approximation and Madhava's correction method.    	

\section{Proof for the $\pi$-series based on Madhava's ideas:}
\label{sec:madhava}
In this section we present Madhava's derivation of the $\pi$-series as presented in \emph{yuktibhasha} following the thread of the proof described in \cite[Section~6.3]{Jyesthadeva2008} retaining the same notation as part of Theorem 1. We call its proof as "\emph{upapatti}" in honour of the tradition of proving that existed in ancient Indian Mathematics \cite{Srinivas2005Proofs}. We elaborate the original proof at three
places and present them as three Lemmas:
\begin{enumerate}
\item Lemma 1, \emph{capikaranam, Making the Arc} Lemma.
As part of this we show that in the limit the sum of the approximation of the arc bits is same as the sum of the arc bits. 
\item Lemma 2, \emph{phala parampara, Series of Results} Lemma.
Here we show that for large $n$ we can neglect the reminder of the geometric expansion to get the $\pi$-series.
\item Lemma 3, \emph{samaghata samkalita, Equal Power Sum} Lemma.
As part of this we present a concise version of the proof given in \emph{yuktibhasha} as a proof by induction with explicit handling of the terms to be neglected using big $O$ notation. 
\end{enumerate} 

In order to stay true to the order of presentation in the original the main result i.e Theorem $1$ is presented first followed by the Lemmas which are invoked in the proof of the Theorem.

\textbf{Theorem 1}. The infinite series for $\pi$ i.e the $\pi$-series is $\frac{\pi}{4}=1-\frac{1}{3}+\frac{1}{5}-\frac{1}{7}\cdots$ \\
\\
\textbf{\emph{upapatti}, Proof}:
\begin{figure}
	\centering
	\includegraphics[width=65mm,height=55mm]{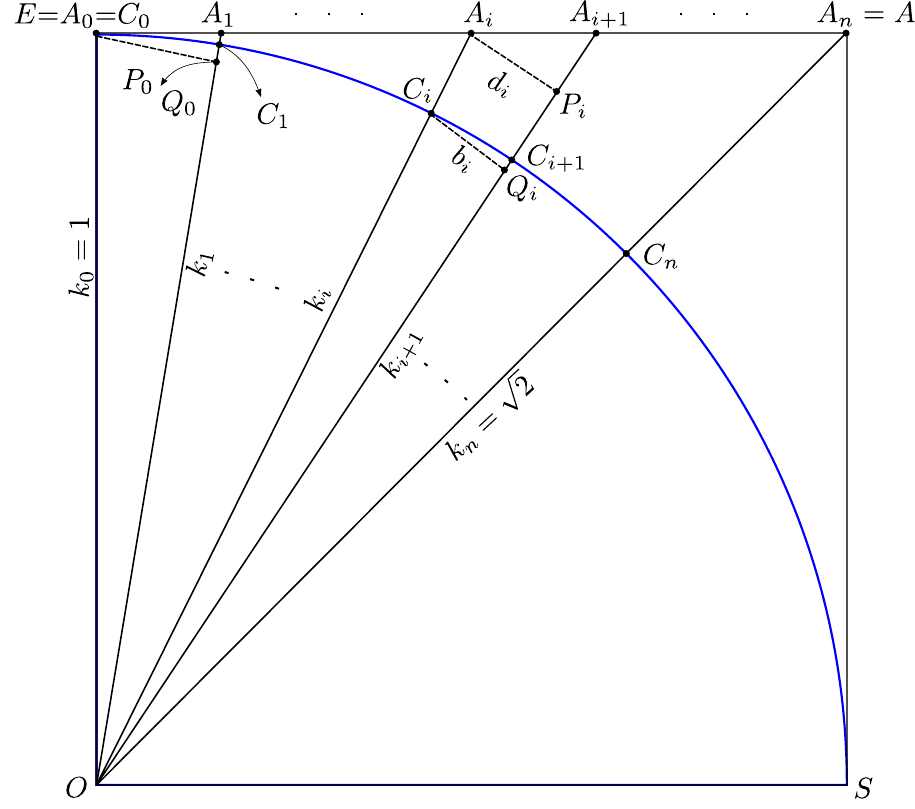}
	\caption{Setting up the arc bit sum for the $\pi$-series. }
	\label{fig:fig_2}
\end{figure}
Consider a circle of radius $1$ with $O$ as its center as shown in Figure \ref{fig:fig_2}. $OEAS$ is a quadrant of the square that circumscribes the circle. 
\begin{itemize}
	\item The tangent $EA$ is divided into $n$ equal parts by marking equidistant points $A_1$, $A_2, \cdots, A_{n}=A$.$ \ \frac{EA}{n}=\frac{1}{n}$. 
	\item The segment $OA_i$ formed by joining the center to $A_i$ is referred to as \emph{karna} i.e a hypotenuse denoted by $OA_i$, of length $k_i$. $i\in\{0,\cdots,n-1\}$. $k_0=1$, $OA=k_n=\sqrt{2}.$
	\item The hypotenuse $OA_i$ (Figure \ref{fig:fig_2}) intersects the circle at $C_i$. 
	$A_iP_i$ is perpendicular to the hypotenuse $OA_{i+1}$. 
	\item  Each $C_iQ_i$ is a \emph{jyardha} or a half-chord or a Rsine obtained by dropping the perpendicular from $C_i$ to $OA_{i+1}$. The points $E$, $C_0$, $A_0$ coincide and so do points $P_0,C_0$ and $C,C_n$.
\end{itemize}

The length of the arc $\overset{\frown}{EC}$ denoted by $\overset{\frown}C$ is $\frac{1}{8}^{th}$ of the circumference
hence $\overset{\frown}C=\frac{\pi}{4}$. 
$\overset{\frown}C$ is the sum of lengths of $C_iC_{i+1}$'s. For large $n$ each $C_iC_{i+1}$ is approximated by the half-chord $C_iQ_i$. 
Let the length of $C_iQ_i$ be denoted by $b_i$ and let the sum of the arc-bit approximations i.e $b_i$'s be denoted by $\hat{A}$: 
\begin{eqnarray*}
\hat{A}&=&b_0 + b_1 + b_2 + \cdots \\
       &=& \lim_{n\rightarrow \infty} \sum_{i=0}^{n-1} b_i    \\
 \end{eqnarray*}  
First an expression for $b_i$ interms of $i$, $n$ and $k_i$ is obtained and then it is shown that for large $n$ adding the arc-bit approximations is same as adding the arc-bits themselves i.e $\overset{\frown}C=\hat{A}$ and finally it is shown that this summation converges to give a series for $\frac{\pi}{4}$.

For $i=0,\cdots,n-1$, triangles $OEA_{i+1}$ and $A_iP_iA_{i+1}$ are similar since they share an angle i.e $EA_{i+1}O$ and $A_iA_{i+1}P_i$ and contain a right angle each i.e $OEA_{i+1}$ and $A_iP_iA_{i+1}$ respectively. Hence we have:
\begin{eqnarray*}
	A_{i}P_i&=& A_iA_{i+1} \frac{OE}{OA_{i+1}}  = \left(\frac{1}{n}\right) \left(\frac{1}{k_{i+1}}\right) \\
\end{eqnarray*} 
To get an expression for $C_iQ_i$ $(b_i)$, similarity of 
triangles $OC_iQ_i$ \&  and $OA_{i}P_{i}$ is considered and the
expression for $A_iP_i$ from the previous step is substituted to get:
\begin{eqnarray*}
	b_i=C_iQ_i&=&A_iP_i \frac{OC_i}{OA_i} \\
	&=& \left(\frac{1}{n}\right)\left( \frac{1}{k_{i+1}}\right) \left( \frac{1}{k_i}\right) \\
	&=&\left(\frac{1}{n}\right)\left( \frac{1}{k_{i}k_{i+1}}\right)
\end{eqnarray*}

Applying \emph{capikaranam, Making the
Arc} Lemma i.e Lemma 1, with the expression for $C_iQ_i$ derived above, we have the arc-bit sum i.e $\overset{\frown}{C}$ for large $n$ equals $\hat{A}$. 
\begin{eqnarray*}
	\overset{\frown}{C} = \hat{A} = \lim_{n\rightarrow \infty} \sum_{i=0}^{n-1} b_i \\
\end{eqnarray*} 
It is shown that the sum to be evaluated does not change if we substitute $k^2_i$ for $k_ik_{i+1}$ as follows:
\begin{eqnarray*}
k_i&=& 1 + \left(\frac{i}{n}\right)^2 \\
k^2_i&\leq& k_ik_{i+1} \leq k^2_{i+1} \\
 \left(\frac{1}{n}\right) \sum_{i=0}^{n-1}\left( \frac{1}{k^2_{i+1}}\right) &\leq&
 \left(\frac{1}{n}\right) \sum_{i=0}^{n-1}\left( \frac{1}{k_{i}k_{i+1}}\right) \leq
 \left(\frac{1}{n}\right) \sum_{i=0}^{n-1}\left( \frac{1^2}{k^2_{i}}\right)\\
\end{eqnarray*}
In the inequality given above the sum on the left and the one on the right are lower and upper-bounds respectively for the sum we are interested in and their difference goes to zero for large $n$ as:
 \begin{eqnarray*}
  \left(\frac{1}{n}\right) \left[  \left(\frac{1}{k^2_n}\right) - \left( \frac{1}{k^2_0}\right) \right] \\
  k^2_0=1, \ k^2_n=2 \\
  \left|\left(\frac{1}{n}\right) \left[  \left(\frac{1}{k^2_n}\right) - \left( \frac{1}{k^2_0}\right) \right] \right| = \frac{1}{2n} \\
 \end{eqnarray*}
Thus if the limit $\lim\limits_{n \rightarrow \infty} \left(\frac{1}{n}\right) \sum_{i=0}^{n-1}\left( \frac{1^2}{k^2_{i}}\right)$ exists it will be same as
$\lim\limits_{n \rightarrow \infty} \left(\frac{1}{n}\right) \sum_{i=0}^{n-1}\left( \frac{1^2}{k^2_{i+1}}\right)$. By
 \cite[Theorem~3.25]{rudin1964principles}, since $\overset{\frown}{C}$ is `sandwiched' between the two sums which go to the same value we can use either one of them to evaluate it. We use
$\left(\frac{1}{n}\right) \sum_{i=0}^{n-1}\left( \frac{1^2}{k^2_{i}}\right)$. Substituting the value of $k_i^2$ in the series gives:
\begin{eqnarray*}
\overset{\frown}{C} &=& \lim_{n\rightarrow \infty} \left(\frac{1}{n}\right) \sum_{i=0}^{n-1} \frac{1}{\left( 1 + \left( \frac{i}{n}\right)^2  \right)} \\
\end{eqnarray*}	
 By \emph{phala parampara,  Series of Results} Lemma i.e Lemma 2, $\overset{\frown}{C}$ is: 
\begin{eqnarray*}
	\overset{\frown}{C}&=&\sum_{p}^{\infty} (-1)^p \frac{1}{2p+1} \\
	&=& \pi-series \\
\end{eqnarray*}
Since the L.H.S is $\frac{1}{8}^{th}$ of the circumference of a circle of unit radius we have:
\begin{eqnarray*} 
	\frac{\pi}{4}&=&\sum_{p}^{\infty} (-1)^p \frac{1}{2p+1} \\
				 &=&1-\frac{1}{3}+\frac{1}{5}-\frac{1}{7}+\cdots \\
\end{eqnarray*}	
Hence the proof.

In the following Lemma we use the trigonometric inequality i.e $tan(x)\ge x$ to show that the sum of the arc bits i.e $\overset{\frown}C$ is same as $\hat{A}$ the sum of the
arc-bit approximations.

\textbf{Lemma 1}. \textbf{\emph{capikaranam, Making the Arc}}.
\label{lem:1}
Given the arc-bits $c_i$, $i\in\{{0,\cdots,n-1}\}$ that denote the length of $C_iC_{i+1}$ and their half-chord approximations i.e $b_i$,
 let the sum of the approximations be denoted by
$\hat{A}=\lim\limits_{n\rightarrow \infty} \sum_{i=0}^{n-1}b_i$ and the sum of the
arc-bits be denoted by $\overset{\frown}{C}=\lim\limits_{n\rightarrow \infty} \sum_{i=0}^{n-1}c_i$, then $\overset{\frown}{C} = 
\hat{A}$.  \\  \\
\textbf{Proof:}   
This Lemma serves as a link between the geometric intuition and the algebra to find the series sum that follows to establish convergence.
We show that the difference between the infinte sums of the
arc-bits and their approximations goes to zero using $c_i=\theta_i$=
$arctan(\frac{bi}{OQ_i}) \leq \frac{bi}{OQ_i}$\footnote{This is due the
trigonometric inequality $tan(x)\ge x, 0< x < \frac{\pi}{2}$, please refer Section \ref{subsub:tanx}
where it is shown that this can be proved without resorting to the area of the sector of a circle.}
Let:
\begin{eqnarray*}
d(n) = \sum_{i=0}^{n-1} \left ( c_i - b_i \right) \\
d(n) = \sum_{i=0}^{n-1} \left ( \theta_i - b_i \right) \\
d(n) = \sum_{i=0}^{n-1} \left ( arctan\left(\frac{b_i}{OQ_i}\right) - b_i  \right) \\
\end{eqnarray*}
Given:
\begin{eqnarray*}
b_i &=& C_iQ_i=\frac{1}{n} \frac{1}{k_ik_{i+1}} \ , \ OQ_i=\sqrt{1-b_i^2} \\
Let \ \lambda_i &=&\frac{1}{k_ik_{i+1}} \\
\frac{b_i}{OQ_i} &=& \frac { \frac{1}{n} \lambda_i}{\sqrt(1-\frac{1}{n^2}\lambda_{i}^2)} \\
arctan(x) &\le& x \ , \ 0\le x \le 1 \\
arctan(\frac{b_i}{OQ_i}) &\le& \frac{b_i}{OQ_i} \\
\end{eqnarray*}
Using the above we have:
\begin{eqnarray*}
	 d(n) &\leq&  \sum_{i=0}^{n-1} \left( \frac { \frac{1}{n} \lambda_i}{\sqrt(1-\frac{1}{n^2}\lambda_{i}^2)} - \frac{1}{n}\lambda_i \right) \\
\end{eqnarray*}
Considering modulus of $d(n)$ and applying the Triangle inequality for modulus of the sum we have: 
\begin{eqnarray*}
 |d(n)| &\leq&  \sum_{i=0}^{n-1} \left| \left( \frac { \frac{1}{n} \lambda_i}{\sqrt(1-\frac{1}{n^2}\lambda_{i}^2)} - \frac{1}{n}\lambda_i \right) \right | \\
 &=&  \sum_{i=0}^{n-1} \left | \frac{  \lambda_i}{n}\left(\frac{1}{\sqrt{1-\frac{1}{n^2}\lambda_i^2}} - 1 \right) \right | \\ 
\end{eqnarray*}
Using:
\begin{eqnarray*}
1 &\leq& k_ik_{i+1}  \leq 2 \\ 
0.5=\lambda_{min} &\leq&  \lambda_i \leq \lambda_{max}=1 	\\
\end{eqnarray*}
Since $\frac{1}{\sqrt{1-\frac{1}{n^2}\lambda_i^2}}$ is $\ge 1$,
we have $ \left | \frac{\lambda_i}{n}\left(\frac{1}{\sqrt{1-\frac{1}{n^2}\lambda_i^2}} - 1 \right) \right |$ is $\ge$ $0$. Since $|d(n)|$ increasees with $\lambda_i$ we can use $\lambda_{max}$ instead of $\lambda_i$ to do away with the summation as:
\begin{eqnarray*}
 |d(n)| &\leq&  n  \frac{\lambda_{max}}{n} \left(\frac{1}{\sqrt{1-\frac{1}{n^2}\lambda_{max}^2}} - 1 \right)\\
\end{eqnarray*}
Taking $\lim\limits_{n\rightarrow \infty} |d(n)|$ we have:
\begin{eqnarray*}
\lim_{n \rightarrow \infty} |d(n)|&=& \lim_{n\rightarrow \infty}  \left(\frac{1}{\sqrt{1-\frac{1}{n^2}}} - 1 \right) \\
&=& 0 \\
\end{eqnarray*}
Thus $\lim\limits_{n\rightarrow \infty} d(n)=0$ which is used to get the desired result as follows:
\begin{eqnarray*}
\overset{\frown}{C}&=&\lim_{n\rightarrow \infty} \sum_{i=0}^{n-1}c_i \\
&=&\lim_{n\rightarrow \infty} \sum_{i=0}^{n-1}(c_i-b_i)+b_i \\
&=&\lim_{n\rightarrow \infty} \sum_{i=0}^{n-1}(c_i-b_i)+\lim_{n\rightarrow \infty}\sum_{i=0}^{n-1}b_i\\
&=& \lim_{n\rightarrow \infty} d(n) + \hat{A} \\
&=& \hat{A}\\
\end{eqnarray*}
We have proved that if the sum of the approximations of the arc-bits converges then it is equal to the sum of the arc-bits themselves. \\

In the Lemma that follows i.e the \emph{phala parampara}, Series of Results Lemma we show that as indicated in
\emph{yuktibhasha} for large $n$ we can stop expanding the geometric series when the reminder becomes negligeble enough to ensure convergence of the series.

\textbf{Lemma 2}, \textbf{\emph{phala parampara, Series of Results}}, The sum of arc-bit approximations
denoted by $\hat{A}$ is same as the 
${\pi}$-series sum, where $\hat{A}$ and $\pi$-series are as defined below:
\begin{eqnarray*}
\hat{A} &=&  \lim_{n\rightarrow \infty} \left(\frac{1}{n}\right) \sum_{i=0}^{n-1}\frac{1}{\left( 1 + \left( \frac{i}{n}\right)^2  \right)} \\
\pi-series&=&\sum_{p=0}^{\infty}(-1)^p \frac{1}{2p+1} \\
\end{eqnarray*}
\textbf{Proof:} \\
Let: 
\begin{eqnarray*}
t(p,n)&=& \frac{1}{n}\sum_{i=0}^{n-1} \left(\frac{i}{n}\right)^{2p} \\
r(p,n)&=& \frac{1}{n}\sum_{i=0}^{n-1} \frac{\left(\frac{i}{n}\right)^{2p}}{1+(\frac{i}{n})^2}\\
T(p)&=& \lim_{n\rightarrow \infty}t(p,n) \\
\end{eqnarray*}
Since $1+\left(\frac{i}{n}\right)^2 \ge 1$ :
\begin{eqnarray*}
r(p,n)&=& \frac{1}{n}\sum_{i=0}^{n-1} \frac{\left(\frac{i}{n}\right)^{2p}}{1+(\frac{i}{n})^2} \le \frac{1}{n}\sum_{i=0}^{n-1} \left(\frac{i}{n}\right)^{2p} \\
r(p,n) &\le& t(p,n) \\
\end{eqnarray*}
By Lemma 3 i.e the \emph{samaghata samkalita, Equal Power Sum} Lemma, which gives the equal-power series sum for large $n$ we have:
\begin{eqnarray*} 
T(p)&=& \lim_{n\rightarrow \infty} \frac{1}{n}\sum_{i=0}^{n-1} \left(\frac{i}{n}\right)^{2p}\\
&=& \frac{1}{2p+1} \\
\end{eqnarray*}
The $\pi$-series can be written as the limit:
\begin{eqnarray*}
	\pi-series=\lim_{m\rightarrow \infty} \sum_{p=0}^{m} (-1)^p \frac{1}{2p+1} = \lim_{m\rightarrow \infty} \sum_{p=0}^{m} (-1)^p T(p)  \\
\end{eqnarray*}
This series sum is well defined since it is an alternating series in which the terms of the series decrease to zero \cite[Theorem~3.43]{rudin1964principles}.
To prove this Lemma we have to show that $\forall$ $\epsilon > 0$, $\exists$ a $M$
s.t $\forall$ $m > M$ the following is true:
\begin{eqnarray*}
	\left| \hat{A} - \sum_{p=0}^{m} (-1)^p T(p)  \right | < \epsilon \\
\end{eqnarray*}
This can be accomplished by writing the geometric series expansion of the term  $\frac{1}{\left( 1 + \left( \frac{i}{n}\right)^2  \right)}$ in $\hat{A}$ as the sum of first $M$ terms and the second part i.e the reminder of the geometric series expansion. 
In \emph{yuktibhasha} the geometric expansion\footnotemark is written explicitly describing the last term with an instruction to stop expanding when the last term i.e the reminder is small enough\cite[Section~6.3.3]{Jyesthadeva2008}. Our proof is based on this instruction.
\footnotetext{
	Let, $x = \left(\frac{i}{n}\right)^2,  i \in \{ 0,\cdots,n-1\}$
\begin{eqnarray*}
	 \frac{1}{1+x}&=&1-x\frac{1}{1+x}=1-x\left(1-\frac{1}{1+x}\right)=1-x+x^2\frac{1}{1+x}=\cdots \\ 
	 \frac{1}{1+x}&=&\sum_{p}^{M-1}(-1)^p x^{2p} + (-1)^M \frac{x^{2M}}{1+x} \\
\end{eqnarray*}
}
\begin{eqnarray*}
\frac{1}{\left( 1 + \left( \frac{i}{n}\right)^2  \right)} &=&
\sum_{p=0}^{M-1} (-1)^p \left(\frac{i}{n}\right)^{2p} + (-1)^M\frac{\left(\frac{i}{n}\right)^{2M}}{1+\left(\frac{i}{n}\right)^2} \\
\end{eqnarray*}
By using the above $\hat{A}$ can be simplified as follows:
\begin{eqnarray*}
\hat{A} &=&  \lim_{n\rightarrow\infty} \left(\frac{1}{n}\right) \sum_{i=0}^{n-1} \sum_{p=0}^{M-1} (-1)^p \left(\frac{i}{n}\right)^{2p} +  (-1)^M \lim_{n\rightarrow\infty} \left(\frac{1}{n}\right) \sum_{i=0}^{n-1} \frac{\left(\frac{i}{n}\right)^{2M}}{1+\left(\frac{i}{n}\right)^2} \\
\end{eqnarray*}
After interchanging the finite limits and substituting for the terms we have:
\begin{eqnarray*}
\hat{A}&=&   \sum_{p=0}^{M-1} (-1)^p  \lim_{n\rightarrow\infty} \left(\frac{1}{n}\right) \sum_{i=0}^{n-1} \left(\frac{i}{n}\right)^{2p}  + (-1)^M \lim_{n\rightarrow \infty}r(M,n) \\
&=&\sum_{p=0}^{M-1} (-1)^p T(p) + (-1)^M \lim_{n\rightarrow \infty}  r(M,n) \\	
\end{eqnarray*}
Comparing $r(M,n)$ and $t(M,n)$ , since $r(p,n) \le t(p,n)$ we have:
\begin{eqnarray*}
 lim_{n\rightarrow \infty}r(M,n) \le lim_{n\rightarrow \infty} t(M,n)	\\
\end{eqnarray*}
Since $lim_{n\rightarrow \infty} t(M,n)=T(M)$ the above proves existence of the L.H.S limit \cite[Theorem~3.25]{rudin1964principles}. Using the above inequality we find the $M$ for a given $\epsilon$ as follows:
\begin{eqnarray*}
	 \hat{A}-\sum_{p=0}^{m} (-1)^p T(p)  &=&
	\hat{A}- \sum_{p=0}^{M-1} (-1)^p T(p) - \sum_{p=M}^{m} (-1)^p T(p)  \\
	&=&\sum_{p=0}^{M-1} (-1)^p T(p) + (-1)^M \lim_{n\rightarrow \infty}  r(M,n) - \sum_{p=0}^{M-1} (-1)^p T(p) - \sum_{p=M}^{m} (-1)^p T(p) 
	\\
	&=& (-1)^M \lim_{n\rightarrow \infty}  r(M,n) - \sum_{p=M}^{m} (-1)^p T(p)   \\
\end{eqnarray*}
 By using the Triangle inequality on the magnitude of
 the difference and the fact that $lim_{n\rightarrow \infty}r(M,n) \le lim_{n\rightarrow \infty} t(M,n)$ we have:
\begin{eqnarray*}
	\left| (-1)^M \lim_{n\rightarrow \infty}  r(M,n)-\sum_{p=M}^{m} (-1)^p T(p) \right| &\le&  \left|\lim_{n\rightarrow \infty} t(M,n)\right| +\left| \sum_{p=M}^{m} (-1)^p T(p) \right| \\
	&\le& \left| T(M) \right| + \left|\sum_{p=M}^{m} (-1)^p T(p)\right|  \\
\end{eqnarray*}
For an alternating series since the sum $\left|\sum_{p=M}^{m} (-1)^p T(p)\right|$ is $\le$ $|T(M)|$ we have:
\begin{eqnarray*}
   \left| T(M) \right| + \left|\sum_{p=M}^{m} (-1)^p T(p)\right| &\le& \left| T(M) \right| + \left|T(M)\right|\\
 &\le& 2\frac{1}{2M+1} \\	
\end{eqnarray*}
We can choose $M$ s.t $\frac{1}{2M+1} < \frac{\epsilon}{2}$ so that $\hat{A}$ is arbitrarily
close to $\pi$-series sum. Thus we have
proved that the sum of the arc-bit approximations denoted by $\hat{A}$ is the ${\pi}$-series.

Lemma 3 is used to approximate the sum of powers of natural numbers for large $n$. The proof for this is given in detail in \emph{yuktibhasha} \cite[Section~6.4]{Jyesthadeva2008}. We present a concise version of this proof as a `proof by induction', writing the error terms explicitly using the big $O$ notation to keep track of the terms that can be neglected for large $n$.

\textbf{Lemma 3}. \textbf{\emph{samaghata samkalita, Equal Power Sum }}. T.S.T $\lim\limits_{n\rightarrow \infty} \frac{1}{n} \sum_{i=0}^{n-1}\left(\frac{i}{n}\right)^{2p} = \frac{1}{2p+1}$\\ \\
\textbf{Proof:}
   
In the notation that we have used in Lemma 2 the limit given in the statement above is:
\begin{eqnarray*}
\lim_{n\rightarrow \infty} t(p,n) &=& T(p) \\
t(p,n)&=& \frac{1}{n} \sum_{i=0}^{n-1}\left(\frac{i}{n}\right)^{2p} \\
T(p) &=& \frac{1}{2p+1}  \\
\end{eqnarray*}
$\lim_{n\rightarrow \infty} t(p,n)$ can be written as:
\begin{eqnarray*} 
\lim_{n\rightarrow \infty} t(p,n) =
\lim_{n\rightarrow \infty} \frac{1}{n^{2p+1}} \sum_{i=0}^{n-1}\left(i\right)^{2p} \\
\end{eqnarray*}
Thus we need to look at the sum of the powers of the first $n$ natural numbers, i.e $S^{p}_n=\sum_{i=1}^{n} i^p$ for large $n$, and then use the result with the power term being an even number ($2p$). 

Statement of Induction: Let $S^p_n$ = $\sum_{i=1}^{n} i^{p}$
\begin{eqnarray*}
	S_n^{p-1} &=& \frac{n^p}{p} + O(n^{p-1})\\	
	\implies S_n^{p} &=& \frac{n^{p+1}}{p+1} + O(n^{p})\\
\end{eqnarray*}
Proof: \\
\begin{eqnarray*}
S^{p-1}_n &=& n^{p-1} + (n-1)^{p-1} + \cdots + 1 \\
S^{p}_n &=& n^{p} + (n-1)^{p} + \cdots + 1 \\
nS_n^{p-1} &=&  n^{p} + n(n-1)^{p-1} + \cdots + n \\
nS_n^{p-1} - S^{p}_n &=& n^p-n^p + [n-(n-1)](n-1)^{p-1} + [n-(n-2)](n-2)^{p-1} + \cdots + (n-1) \\ 
nS_n^{p-1} - S^{p}_n &=& (n-1)^{p-1} + 2(n-2)^{p-1}  + 3 (n-3)^{p-1} +\cdots + (n-1) \\
\end{eqnarray*}
The idea is to write $nS_n^{p-1} - S^{p}_n$ in terms of sums of powers upto
$n-1$, $n-2$ $\cdots$, $1$  and then use $S_n^{p-1}$ given in the premise of the statement to get an expression for $S^p_n$ for large $n$. 
$nS_n^{p-1} - S^{p}_n$ can be expanded as a telescopic sum as follows:

\begin{eqnarray*}
	nS_n^{p-1} - S^{p}_n = (n-1)^{p-1} + (n-2)^{p-1} +  &+& (n-3)^{p-1} +\cdots + 1 \\
	 						 	 + (n-2)^{p-1} &+&  (n-3)^{p-1} + \cdots + 1 \\
	 						      			   &+& (n-3)^{p-1} + \cdots + 1 \\	
	 						      			   &\vdots& \\
	 						      			   &\cdots& \ \cdots \ \cdots \ \cdots \ \cdots  \ \ \ \ + 1 \\
	 nS^{p-1}_n - S^{p}_n = 	S^{p-1}_{n-1} + S^{p-1}_{n-2} + \cdots + S_{n-(n-1)}^{p-1}  \\	      			   
\end{eqnarray*}
From the statement of Induction using $S_n^{p-1} = \frac{n^p}{p} + O(n^{p-1})$
\begin{eqnarray*}
nS^{p-1} - S^{p}_n &=& \frac{(n-1)^p}{p} + O(n^{p-1}) + \frac{(n-2)^p}{p} + O(n^{p-1})
+\cdots + \frac{(n-(n-1))^p}{p} + O(n^{p-1}) \\
nS_n^{p-1} - S^{p}_n &=& \frac{(n-1)^p}{p}  + \frac{(n-2)^p}{p} 
+\cdots + \frac{(n-(n-1))^p}{p} + nO(n^{p-1}) \\
&=& \frac{S^{p}_{n-1}}{p} + \frac{n^p}{p} - \frac{n^p}{p} + O(n^{p})  \\
&=& \frac{S^{p}_{n}}{p} + O(n^{p}) \\  
\end{eqnarray*}
We have, $nS_n^{p-1} - S^{p}_n = \frac{S^{p}_n}{p} + O(n^p)$.
Solving for $S^{p}_n$:
\begin{eqnarray*}
\left( \frac{p+1}{p} \right)S^{p}_n &=& n S_n^{p-1} + O(n^p) \\
\left( \frac{p+1}{p}\right)S^{p}_n &=& n \frac{n^p}{p} + O(n^p) \\
S^{p}_n &=& \frac{n^{p+1}}{p+1} + O(n^p)
\end{eqnarray*}
Thus for large $n$ we have the approximation of $S^{p}_n \approx \frac{n^{p+1}}{p+1} $

To set the induction into motion we show that $S^1_n$ = $\frac{n^2}{2}$ + $O(n)$
This is can be shown by using the expression for sum of first $n$ natural numbers which was known to ancient Indian mathematicians like Aryabhata \cite[Section~6.4.5]{Jyesthadeva2008} as: $\frac{n(n+1)}{2}=\frac{n^2}{2}+\frac{n}{2}$, that gives $S^1_n=\frac{n^2}{2}+O(n)$.
Now we connect this result to the statement of this Lemma. We need to show that $\lim_{n\rightarrow \infty} t(p,n)=T(p)=\frac{1}{2p+1}$. $t(p,n)$ involves summing
upto $n-1$ we show that for large $n$ we can use the
same expression for $S^p_{n-1}$ as that of $S_n^{p}$.
\begin{eqnarray*}
	S^p_{n-1} &=&  S^{p}_n - n^p \\
			 &=&  \frac{n^{p+1}}{p+1} + O(n^p) - n^p \\
			 &=&  \frac{n^{p+1}}{p+1} + O(n^p)=S^p_n \\
\end{eqnarray*}
Using this in $t(p,n)$ we have:
\begin{eqnarray*}
	t(p,n)&=& \frac{1}{n^{2p+1}}\left(\frac{n^{2p+1}}{2p+1} + O(n^{2p})\right)  \\
	&=& \frac{1}{2p+1} + \frac{O(n^{2p})}{n^{2p+1}} \\
\end{eqnarray*}
For large $n$, $\frac{O(n^{2p})}{n^{2p+1}}$ which is $\le \frac{K}{n}$ ($K$ a constant $> 0$) can be made as small as needed.
\begin{eqnarray*}
	\left|t(p,n)-\frac{1}{2p+1}  \right| &=& \frac{O(n^{2p})}{n^{2p+1}} \\
	&\le& \frac{K}{n}  \\  
	\lim_{n\rightarrow \infty} t(p,n) &=& T(p) = \frac{1}{2p+1}   \\
\end{eqnarray*}
Hence the proof.

\section{The arctan, sine and cosine series.}
\label{sec:madhava_other}
Along with the $\pi$-series \emph{yuktibhasha} also presents derivations of the arctan, sine and cosine series. We show that the fundamental method used in all these derivations is the same, we briefly describe the steps involved in deriving the sine series and compare them with the corresponding steps in the derivation of the $\pi$-series. 
\subsection{The arctan series}
\label{subsec:madhava_arctan}
The derivation of the $\pi$-series starts by dividing $EA=1$ (Figure \ref{fig:fig_2}) into
$n$ equal parts. Instead if we start by dividing $EA=x, 0\le x \le 1$ into $n$ equal parts we get the the $\arctan$ series. 
The derivation of the series for $\pi$ can be viewed as a specific case of the derivation when the length of $EA$ is $x$ where $0 \le x \le 1$,  
giving the $\pi$-series when $x=1$.
With $C$ as the point where $OA$ intersects the circle, let $\theta$
represent the arc $\overset{\frown}{EC}$, then $tan({\theta})=\frac{EA}{1}=x, \theta= \arctan(x)$ \cite[Section~6.6]{Jyesthadeva2008}. If we divide 
$EA$ into $n$ equal parts each of length $\frac{x}{n}$, $\frac{x}{n}$ figures
in the derivation of Theorem 1 instead of $\frac{1}{n}$ giving:
\begin{eqnarray*}
\overset{\frown}{EC}=\theta=arctan(x)=x-\frac{x^3}{3}+\frac{x^{5}}{5}+\cdots
\end{eqnarray*} 
\subsection{The sine and cosine series.}
\label{subsec:madhava_sine}
\begin{figure}
	\centering
	\includegraphics[width=60mm,height=50mm]{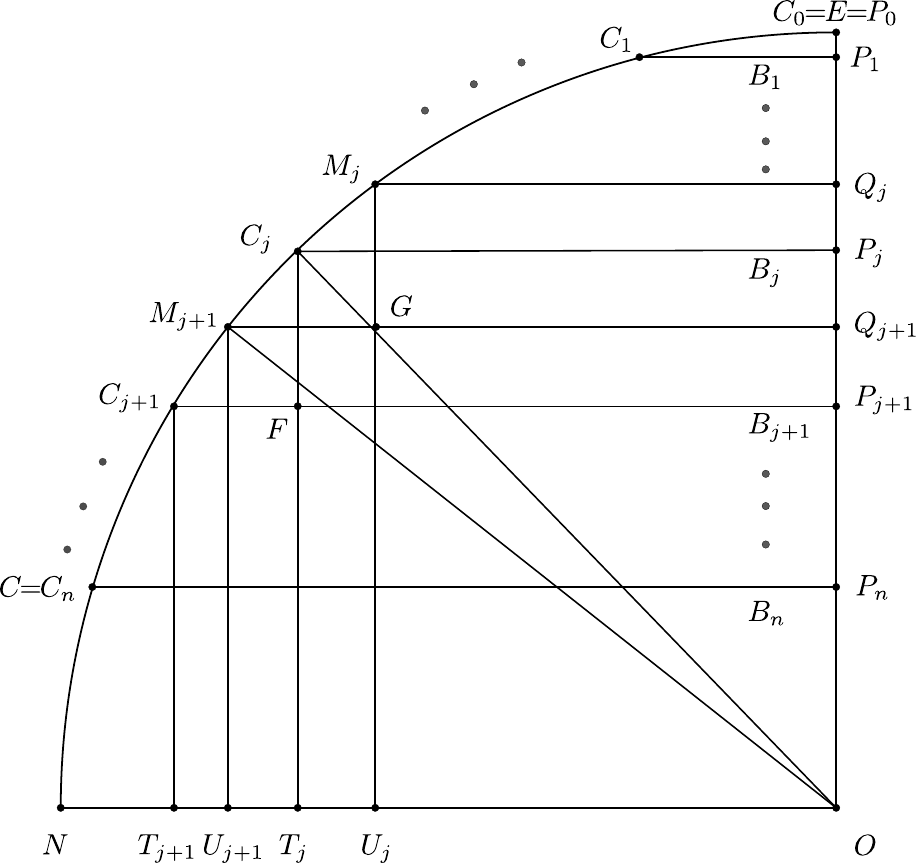}
	\caption{The sine series \emph{bhuja} bits}
	\label{fig:fig_3}
\end{figure}
\emph{Yuktibhasha} describes Madhava's derivation of the sine and the versine i.e 1-cosine
series as well. Without going into the details of the proof we show that the
method used for deriving these series is fundamentally similar to that of the
$\pi$-series. The method that Madhava followed involves the following steps:
\begin{enumerate}
\item Writing the quantity of interest like the arc of a circle, sine or versine as sum of $n$ bits, where the bits become small as $n$ increases.
\item Setting up an iterative process of refining the estimate of these 
bits starting with an initial estimate which gives the series.  
\item  Approximating the sum of powers of natural numbers that occur as part of the series for large $n$. 
\item The process of refinement and approximation of the sums leads to a monotone decreasing alternating series so that the reminder is always in sight and can be seen to be decreasing as the iterations proceed. 
\end{enumerate}
We briefly explain the steps in the derivation of the sine series ( Figure \ref{fig:fig_3} )  \cite[Section~7.4.3]{Jyesthadeva2008}, highlighting the differences as compared to the $\pi$-series. The differences are 
tabulated in Table \ref{tab:sine-comparision}.

The arc $\overset{\frown}{EC}$ of length $s$ is divided into $n$ equal parts each of length $\frac{s}{n}$ given by $EC_1, C_1C_2,\cdots, C_jC_{j+1}\cdots,C_{n-1}C$. The `bhuja-bit' $B_j=C_jP_j$ is   $sin(\frac{js}{n})$. We are interested in finding $B_n=sin(s)$. $B_n$ is written as the sum of successive differences of sines, i.e $B_j$'s as:
\begin{eqnarray*}
	 \Delta_j&=&B_j-B_{j-1} \\
	 B_n&=&\sum_{j=1}^{n} \Delta_j
\end{eqnarray*}
Thus we have the first step in Madhava's method i.e the sum of $n$ bits which shrink as $n$ increases. To set up the iteration an important result pertaining to
second order difference of sines attributed to Arybhata\footnote{Bhaskara used first order sine differences i.e $sin(y')-sin(y)=(y'-y)cos(y)$ akin to differentiation to calculate instantaneous planetary motion \cite[Chapter~8]{srinivasiengar1967history}} \cite[Section~5.3]{Ramasubramanian2010} i.e $\Delta_j-\Delta_{j-1}=(\alpha)^2B_j, \alpha=2 \  sine(\frac{s}{2n})$ is used to show that $B_n=B(s)$ can be
expressed as:
\begin{eqnarray*}
	B_n = nB_1- \left( \alpha \right)^2 [(B_1+\cdots+B_{n-1}) + (B_1+\cdots+B_{n-2})+\cdots+1] \\
\end{eqnarray*}
\begin{table}[htbp]
	\centering
	\caption{Comparision of the steps in Madhava's method for $\pi$ and 
		sine series}
	\label{tab:sine-comparision}
	\begin{tabular}{lcc} 
		\toprule 
		\textbf{Steps in Madhava's method} & \textbf{$\pi$-series} &
		\textbf{sine series} \\
		\toprule
		Sum for large $n$ &     	$\hat{A}=\sum_{i=0}^{n-1}b_i$       & $B=\sum_{i=0}^{n-1}\Delta_i$ \\
		&	$b_i$ is the arc-bit approximation.&
		$\Delta_i$ is the difference of sines \\
		& & \\
		Iterative refinement & 	By Geometric Series Expansion& 	By Second Order Difference of Sines \\
		& & \\
		Power Series in $m^{th}$ iteration  &    $S^{2m}_n \approx \frac{n^{2m+1}}{2m+1}$ & 
		$S^{2m-1}_{n-1}+\cdots+S^{2m-1}_{1}\approx\frac{n^{2m+1}}{2m+1!}$                  \\
		& & \\
		The series & 		  			$\frac{\pi}{4}=1-\frac{1}{3}+\frac{1}{5}+\cdots$ & $\sin(\theta)=\theta-\frac{\theta^{3}}{1.2.3}+\frac{\theta^5}{1.2.3.4.5}+\cdots$ \\
		
		\bottomrule 
	\end{tabular}
\end{table}
In the second step of the method the bits i.e $B_j$'s are iteratively refined. 
If the current approximation of sine is
substituted in $B_1,\cdots,B_{n-1}$ i.e in the R.H.S of the above equation we get $B_n$ a refined estimate of sine.
To initiate the refinement $B_j$ 
is assumed to be the arc-bit itself i.e $B_j=\frac{js}{n}$ giving rise to the
following sum:
\begin{eqnarray*}
B_n \approx s- \left(\frac{s}{n}\right)^2 [S_{n-1}+S_{n-2}+\cdots+1]\\
\end{eqnarray*}
After approximating the sum of the powers of natural numbers
for large $n$   
i.e $S^1_{n-1}+S^{1}_{n-2}+\cdots+1$ as $\frac{n^3}{1.2.3}$\footnote{ This sum was known to Narayana Pandita \cite[Section~6.4.5]{Jyesthadeva2008} } we have:
\begin{eqnarray*}
B_n(s) = s-\frac{s^3}{1.2.3}	\\
\end{eqnarray*}
The approximation of the power sum is the third step of Madhava's method.
With an initial estimate of  $B_j=\frac{js}{n}$ we now have a refined estimate of $B_j$ i.e $js-\frac{(\frac{js}{n})^3}{1.2.3}$. There is an addition of a term
because $B_n$ has two terms, the first is $s$, the second multiples the sum's approximation by $s^2$ shifting the powers of $s$ in the previous estimate of sine by $2$. Hence if we substitute this refinement we get a series with three terms and so on. This recursive substitution
plays a role similar to that of the geometric series in the derivation
of the $\pi$-series. The refinement is recursive since the map being estimated using $B_n=sine(s)$ i.e `sine' itself figures in the terms i.e $B_j=sine(\frac{js}{n})$ used to make the refinement. The cosine (versine) series is
discovered similarly by using the refinement of $B_n$.

\section{Proof for the $\pi$-series based on Leibniz's ideas}
\label{sec:leibniz}
In this section we present the popular modern proof for the $\pi$-series which is based on Gottfried Wilhelm Leibniz's ideas. We describe Leibniz's ideas of `Transmutation'and `Integration by parts' using which a rational function is identified to find the area under the arc of a circle \cite{Roy1990Discovery}:
\begin{figure}
	\centering
	\includegraphics[width=60mm,height=35mm]{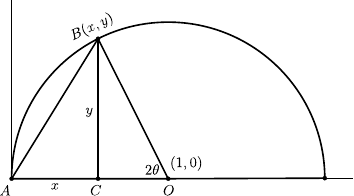}
	\caption{ Area under the arc}
	\label{fig:fig_4}
\end{figure}
Although it was known to Leibniz that a definite integral gives area under a curve, in order to apply it to a circle (Figure \ref{fig:fig_4}) \cite[Section~2]{Roy1990Discovery}, a way to integrate a function like $\sqrt{2x-x^2}$ was not known. Integration of power functions i.e $x^n$ and differentiation of power functions even when the powers were not integers were known \cite[Chapters~18,19]{boyer1968history}. Nicolaus Mercator had just published a way to expand functions like $\frac{1}{1+x}$ as a series of powers of $x$ i.e the geometric series \cite[Section~2]{Roy1990Discovery}. Using all these ingredients Leibniz influenced by ideas of Blaise Pascal developed a method referred to as `Transmutation', a procedure to find a function which could be integrated to get the desired area under the curve. This was applied to the circle, to get a series for $\pi$. 
Leibniz's idea has two parts :
\begin{figure}
	\centering
	\includegraphics[width=80mm,height=70mm]{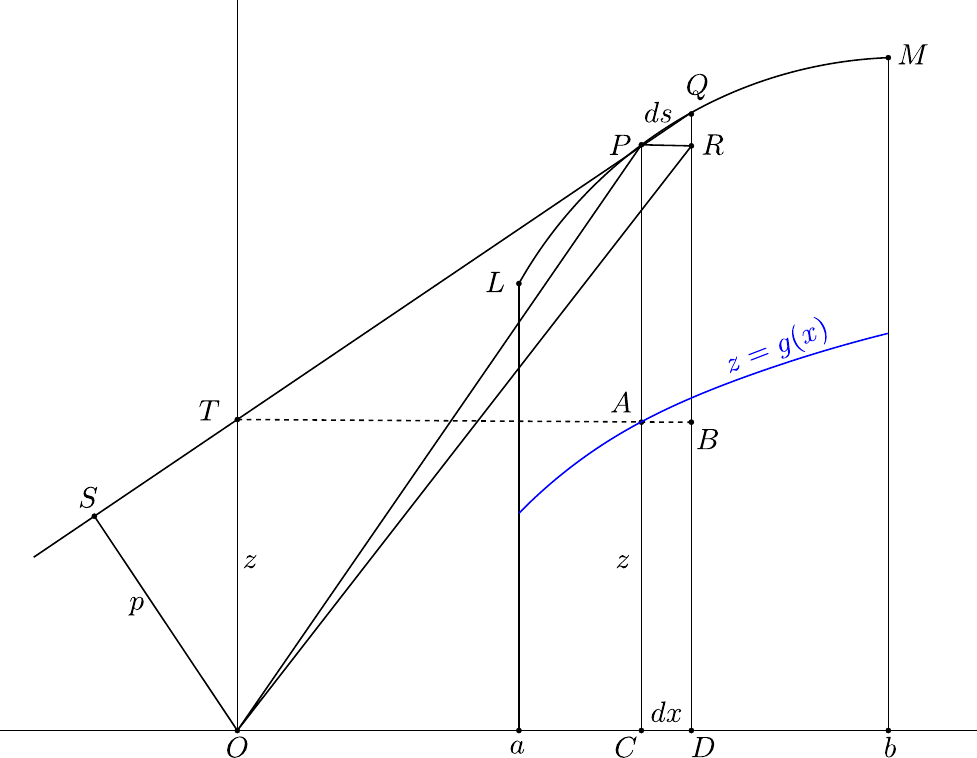}
	\caption{ Leibniz's `Transmutation' of infinitesmal triangluar areas }
	\label{fig:fig_5}
\end{figure}
\begin{enumerate}
\item For the function $y=f(x)$ (Figure \ref{fig:fig_5}) it gives the closed
region bound by the segments $OL$, $OM$ and the the curve $\overset{\frown}{LM}$ as: 
\begin{eqnarray*}
	area \ of \ (O\overset{\frown}{LM}) =\frac{1}{2}\int_{a}^{b} z(x) dx \\
	z(x)=y-x\frac{dy}{dx}
\end{eqnarray*}   
Here $z(x)$ is the
$y$-intercept of the tangent to the function at $x$. 

It was proved as follows:
The characteristic triangle $PQR$ and the triangle $OST$ are similar which gives $\frac{dx}{p}=\frac{ds}{z}$. This means that the area of the triangle $OPQ$ which is $\frac{1}{2}ds \ p$ is equal to $\frac{1}{2}z \ dx$, equals the area of the rectangle $ABCD$. Thus the area of the triangle $O\overset{\frown}{LM}$ visualized as the sum of infinitesmal areas like $OPQ$ is given by $\frac{1}{2}\int_{a}^{b} z(x) dx$. \\ 
\item A result which is similar to integration by parts linking integrals and differentials was derived. Given $y=f(x)$ it is stated as:
\begin{eqnarray*}
	\int_{a}^{b} ydx = [xy]_{a}^{b} - \int_{f(a)}^{f(b)} xdy \\
\end{eqnarray*}
The above result was derived using the first part and by writing the 
area under the curve of $y=f(x)$ i.e under $\overset{\frown}{LM}$ as $\frac{1}{2}f(b)-\frac{1}{2}f(a)+$	area of $(O\overset{\frown}{LM})$. Here $\frac{1}{2}f(b)$ is the area of the triangle $OMb$, $\frac{1}{2}f(a)$ is the area of the triangle $OLa$. As shown above we have, area of ($O\overset{\frown}{LM})=\frac{1}{2}\int_{a}^{b} z(x) dx$. Substituting $z(x)=y-x\frac{dy}{dx}$ for $z(x)$ and evaluating area under ($\overset{\frown}{LM}$) gives the desired rule which is similar to integration by parts.
\end{enumerate}
The two results given above are used to get the  
area of the circle's sector $O\overset{\frown}{AB}$ (Figure \ref{fig:fig_4}) as follows: 

area ($O\overset{\frown}{AB}$) = area of the triangle AOB + area of the closed region between the chord $AB$ and the arc $\overset{\frown}{AB}$ 

The second area i.e the area between the chord and the arc is obtained using Leibniz's transmutation result as $\frac{1}{2}\int_{0}^{x} z dx$ : 
\begin{eqnarray*}
The \ area \ of \ the \ sector \ O\overset{\frown}{AB}= \frac{1}{2}y + \frac{1}{2}\int_{0}^{x} z dx \\
\end{eqnarray*} 
To evaluate the integral we use the second result of Leibniz
with $y=z(x)$ to get: \\
\begin{eqnarray*}
\int_{0}^{x} z dx &=& zx - \int_{0}^{z} x dz \\
\end{eqnarray*}
Using $z=y-x\frac{dy}{dx}$ with $y=\sqrt{2x-x^2}$ we get
$y=z(2-x)$ and $x=\frac{z^2}{1+z^2}$. Substituting these values in the expression for `area of the sector $O\overset{\frown}{AB}$' we get:
\begin{eqnarray*} 
	The \ area \ of \ the \ sector \ O\overset{\frown}{AB} = z - \int_{0}^{z} \frac{t^2}{1+t^2} dt  \\
\end{eqnarray*}
Using the geometric expansion for $\frac{t^2}{1+t^2}$ as $t^2 \sum_p (-1)^p t^{2p}$ in the integral and evaluating area of the sector $O\overset{\frown}{AB}$ as $\theta$ we get the series:
\begin{eqnarray*}
 \theta=z-\frac{z^3}{3}+\frac{z^5}{5}+\cdots 
\end{eqnarray*}
James Gregory, Leibniz's contemporary had also discovered the same series for $arctan$ but had followed a method that is similar to the way coefficients are obtained in Taylor's series, i.e by repeated differentiation. Also he did not substitute $z=1$ in the series, which Leibniz did although strictly speaking he should not have ! since the geometric series does not converge at $z=1$. But his intuition was eventualy proved correct by Abel's result (Section \ref{subsec:pi_abel}). This gave the series for $\frac{\pi}{4}$.
\subsection {The first step, Using the Fundamental Theorem of Calculus}
\label{subsec:ftc}
Leibniz's and Isaac Newton's ideas eventually led to the
the Fundamental Theorem of Calculus \cite[Chapter~1]{boyer1959calculus}, \cite{Bressoud2011} \cite[Theorem~6.20]{rudin1964principles}. While Newton's approach to find a series for $\pi$ was based on the Binomial expansion of
$(1-x)^\frac{1}{2}$ \cite[Section~3]{Roy1990Discovery} the FTC based modern proof encapsulates all that Leibniz attempted using transmutation and integration by parts to find a rational function which could be integrated, giving the first step of the modern day proof for the $\pi$-series: 
\begin{eqnarray*}
	F(x) &=& \int_0^x \frac{1}{1+t^2} dt \ , \  0 \le x \le 1 \\
	F'(x) &=& \frac{1}{1+x^2} \\
	F(x) &=& arctan(x) + C \ \ (F(0)=0 \implies C=0) \\
	arctan(x)&=&\int_{0}^{x}\frac{1}{1+t^2}dt \\
\end{eqnarray*}
\subsection{The series sum, $\pi$-series by analyzing the Reminder after Integration}  
\label{subsec:pi_rem}
The next step in the proof is to expand $\frac{1}{1+t^2}$
 using first $n$ terms of the geometric series and show that the integral of the reminder is negligeble for large $n$.
 \begin{eqnarray*}
 \frac{1}{1+t^2}&=&\sum_{m=0}^{n}(-1)^{m}t^{2m} + (-1)^{n+1} \frac{t^{2n+2}}{1+t^2} \\
 arctan(x)&=&\int_{0}^{x}\frac{1}{1+t^2}dt \ ,\ \ \ \  0 \le x \le 1\\
&=& \int_{0}^{x} \sum_{m=0}^{n} (-1)^m t^{2m} + (-1)^{n+1}\int_0^x \frac{t^{2n+2}}{1+t^2} dt \\
&=&  \sum_{m=0}^{n} (-1)^m \frac{x^{2m+1}}{2m+1} + (-1)^{n+1}\int_0^x \frac{t^{2n+2}}{1+t^2} dt \\
\end{eqnarray*}
Since the magnitude of the integral of the reminder term goes to zero we have the arctan(x) series, with $x=1$ we get the series for $\frac{\pi}{4}$
\cite[Theorem~3.25]{rudin1964principles}:
\begin{eqnarray*}
\left|(-1)^{m+1} \int_0^x \frac{t^{2m+2}}{1+t^2} dt \right| &\le&| \int_0^x t^{2n+2} dt | = \frac{|x^{2n+3}|}{2n+3}\rightarrow 0 \\ 
arctan(x) &=& x -\frac{x^3}{3}+\frac{x^{5}}{5}+\cdots \\
arctan(1) = \frac{\pi}{4} &=&  1 -\frac{1}{3}+\frac{1}{5}+\cdots
\end{eqnarray*}

\subsection{The series sum, $\pi$-series by using Abel's Limit Theorem}
\label{subsec:pi_abel}
Alternately the infinite geometric sum could be considered at
once for $0< x < 1$ when geometric sum converges, the integral can be passed over the sum using the Dominant Convergence Theorem i.e DCT \cite[Theorem~11.32]{rudin1964principles} and for $x=1$ since the $arctan(x)$ is well defined and so is $\frac{1}{1+x^2}$, we can use Abel's Limit theorem \cite[Theorem~8.2]{rudin1964principles} pertaining to power functions to justify that the result holds at $x=1$ as well giving the series for $\pi$ as follows:
\begin{eqnarray*}
arctan(x)&=&\int_{0}^{x}\frac{1}{1+t^2}dt \ ,\ \ \ \  0 < x < 1\\
&=& \int_{0}^{x} \sum_{m=0}^{\infty} (-1)^m t^{2m} dt \\
\end{eqnarray*}
Now to interchange the limits since $|(-1)^m t^{2m}| < 1$ using DCT we have:
\begin{eqnarray*}
arctan(x) &=&	\sum_{m=0}^{\infty} (-1)^m \int_{0}^{x} (-1)^m t^{2m} dt , \ \ 0 < x < 1 \\
&=&	\sum_{m=0}^{\infty} (-1)^m \frac{x^{2m+1}}{2m+1} \\
\end{eqnarray*}
By Abel's Theorem, since a power function sum is involved, since $arctan(x)$ and $\frac{1}{1+x^2}$ are both well defined at $x=1$ the interchange of the limit and the sum remains valid to give the $\pi$-series.
\begin{eqnarray*}
arctan(1)=\frac{\pi}{4}&=& 1 -\frac{1}{3}+\frac{1}{5}+\cdots \\
\end{eqnarray*}

\section{Madhava's and Leibniz's proofs, a Comparision:}
\label{sec:comp}
Now that we have both the proofs i.e the proof based on Madhava's ideas and the popular proof based on
Leibniz's ideas we can compare the proofs. 
We make a step by
step comparision of the proofs highlighting the differences
mainly in terms of the advanced concepts from Calculus used
to derive the proofs, about the role of the trigonometric inequality $tan(x)\ge x$, sum of powers of numbers for large $n$, derivation of other series using the method behind the proofs and Madhava's correction-terms. 
\subsection{$\pi$ as arc length compared to $\pi$ as area
	 under the arc}
\label{subsec:pi_arc}
Madhava's approach adds the arc-bit approximations to 
get to $\pi$ as opposed to Leibniz's approach which is to find the area under the arc. Although $\pi$ figures as a constant in both i.e in the 
circumference and the area of a circle its primary 
association is with the length of the arc as $\pi$ radians
is known to be $180^{\circ}$. This point of departure
at the outset makes the proofs fundamentally different although there are similarities as well in the way the geometric expansion is used and in the way the reminder of the terms 
after summation (Integration) is neglected to prove convergence.  
\subsection {The first step}
\label{subsec:first}
The first step in Leibniz's proof is the main talking point, which is:
\begin{eqnarray*}
arctan(x)=\int_{0}^{x}\frac{1}{1+t^2}dt
\end{eqnarray*}
The equivalent step in Madhava's proof
is :
\begin{eqnarray*}
	\overset{\frown}{C}=\hat{A} &=&  \lim_{n\rightarrow \infty} \left(\frac{1}{n}\right) \sum_{i=0}^{n-1}\frac{1}{\left( 1 + \left( \frac{i}{n}\right)^2  \right)} \\
\end{eqnarray*}
After the inital geometric arguments based on the similarity of triangles the main result that gets us to this stage is 
Lemma 1, where we show that $\overset{\frown}{C}=\hat{A}$. This is based on the trigonometric inequality $tan(x)\ge x$, $0< x < \frac{\pi}{2}$. On the other hand what makes the first step possible in case of Leibniz's proof is the Fundamental Theorem of Calculus (FTC). This difference is significant as FTC is an advanced concept whose derivation in its current form is due to the works of generations of mathematicians mainly, Gregory, Leibniz, Newton, Cauchy, Weierstrass and Reymond \cite[Chapter~7]{boyer1959calculus}. Its derivation is based on the sound knowledge of advanced concepts of Calculus such as Limits, Convergence, Continuity, Differentiation and Integration. Leibniz and Newton laid the foundation, Riemann, Cauchy and Weierstrass rigorously defined the notions of convergence, continuity and limits that is part of modern mathematics \cite[Chapter~7]{boyer1959calculus}. Du Bois Reymond eventually put it in its current form \cite{Bressoud2011}. Out of these concepts Madhava's \emph{upapatti} to get to the first step uses only the concept of limits and convergence. The reason is that it adds infinitesmal lengths where as if we add infinitesmal areas as Leibniz did we need FTC before we start integrating to get the series.

\subsubsection {Role of the trigonometric inequality $\mathbf{tan(x) \ge x \ge sin(x)}$}
\label{subsub:tanx}
In Madhava's proof the role of the inequality is explicit as shown in Lemma 2
where it is used to show that the sum of approximations converges to the sum of arc-bits. Its role is implicit in Leibniz's FTC based first step, where it is used in the derivation of $\frac{d}{dx} arctan(x)$. To be specific it is used in the derivations of differentials of $sin(x)$ and $cos(x)$ which in turn are used in the derivation of the differential of $arctan(x)$\footnote{Please refer \cite[Section~3.3, Section~3.5]{stewart2008calculus} for derivations of differentials of trigonometric functions $sin(x),cos(x),tan(x)$ and $arctan(x)$} .The inequality acts as a crucial bridge between the curve and its approximation in Madhava's proof.
\begin{figure}
	\centering
	\includegraphics[width=35mm,height=70mm]{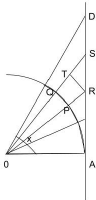}
	\caption{$tan(x)\ge x$ }
	\label{fig:fig_6}
\end{figure}

Briefly let us take a look at what it takes to prove $tan(x)\ge x$. We use the proof given in Stewart's Calculus \cite[Appendix~H]{stewart2008calculus} (Figure \ref{fig:fig_6}) as the outline for the proof. 
It is shown that the $n$ chord-bits whose sum approximates the arc are all less than or equal to the corresponding sub-segments on the tangent that gives $tan(x)$:
\begin{eqnarray*}
	PQ < RT < RS
\end{eqnarray*}
The arc is defined as the sum of these chord-bits ($L_n$) in the limit i.e $x= \sup_n L_n$.

We make one important change ! Instead of defining $x$  as $\sup_n L_n$ we show that $x=\lim\limits_{n \rightarrow \infty} L_n$ using the series expansion for $sine$ which Madhava had derived (Section \ref{subsec:madhava_sine}) and hence was known to Madhava :
\begin{eqnarray*}
L_n &=& 2n \ sine\left(\frac{x}{2n}\right) \\
&=&  2n \left(\frac{x}{2n} - \frac{x^3}{(2n)^33!} + \frac{x^5}{(2n)^55!} \cdots  \right) \\
&=& x + O\left(\frac{1}{n^2}\right) \\
\lim_{n \rightarrow \infty} L_n &=& \lim_{n \rightarrow \infty} x + O\left(\frac{1}{n^2}\right)\\
&=& x 
\end{eqnarray*}
This ensures that  Madhava's proof for the $\pi$-series is not dependent on the proof for area of the sector of a cirlce.

\subsection {The $\pi$-series sum}
\label{subsec:comp_series}
If we use the first option for the series sum in Leibniz's proof (Section \ref{subsec:pi_rem}) which is to show that the integral of the reminder
of the geometric expansion is negligeble then this part is similar to the Lemma 2 i.e the \emph{phala parampara, Series of Results} Lemma in case of Madhava's proof. But there is another popular alternative based on Abel's Limit Theorem (Section \ref{subsec:pi_abel}) that is often used to get the series sum at $x=1$ giving $\frac{\pi}{4}$. In this case additional concepts i.e Dominated Convergence Theorem (DCT) and Abel's result will be needed to complete Leibniz's proof. We tabulate the differences in terms of number
of concepts from advanced Calculus used in the two proofs\footnote{Note that DCT and Abel's Theorem in the table are optional as discussed in Section \ref{subsec:pi_rem}} in Table \ref{tab:comparision}
\begin{table}[htbp]
	\centering
	\caption{Comparision of advanced concepts from Calculus in Madhava's and Leibniz's proofs}
	\label{tab:comparision}
	\begin{tabular}{lcc} 
		\toprule 
		\textbf{Concepts of Mathematical Analysis} & \textbf{Leibniz's Proof} &
		\textbf{Madhava's Proof} \\
		\toprule
		Trigonometric inequality & 				  \checkmark & \checkmark \\
		
		Infinte Series Sum & 			  \checkmark & \checkmark \\
		Convergence Definitions & 		  \checkmark & \checkmark \\
		Geometric Series &				  \checkmark & \checkmark \\
		Alternating Series &			  \checkmark & \checkmark \\
		Continuity &					  \checkmark & $\times$ \\
		Differentiation &				  \checkmark & $\times$ \\
		Integration &					  \checkmark & $\times$ \\
		Fundamental Theorem of Calculus & \checkmark & $\times$ \\
		DCT and Abel's Limit Theorem & \checkmark & $\times$ \\	    				
		\bottomrule 
	\end{tabular}
\end{table}	
\subsection {The sum of Powers of Natural numbers}
\label{subsec:powers}
In case of Leibniz's proof the corresponding step is that of integrating power functions, i.e $x^{2t}$ which was well
known in the $17^{th}$ century. The derivation was
based on differentials of power functions. Where as Madhava used the results of sums of powers known to Indian mathematicians before him for large $n$.
These sums were derived as shown in the \emph{samaghata samkalita, Equal Power Sum} Lemma, i.e Lemma 3 i.e $\sum_{i=1}^n i^p \approx \frac{n^{p+1}}{p+1} $.Thus Madhava
\cite[Section ~6.4.5]{Jyesthadeva2008} knew how to approximate for large $n$ the sum of what is today referred to as the Faulhaber series, named after Faulhaber a German mathematician who worked on it in the $17^{th}$ century. It
was Jacobi and Bernoulli who eventually came up with the coefficients in the $19^{th}$ century \cite{Knuth1993Faulhaber}.

\subsection {Madhava's Correction terms}
\label{subsec:correction}
While Madhava's proof for the $\pi$-series has not received the
attention that it deserves his work on correction terms for the series to hasten its convergence has received considerable attention \cite{hayashi1990correction},\cite{Krishnachandran2024Madhava} mainly because the series is known for converging very slowly \cite{Krishnachandran2024Madhava}. Gregory and Leibniz did not suggest
remedial ways to address this slow convergence \cite[Section~3]{Roy1990Discovery}. 
Madhava arrived at a method to find correction terms based on a unique criterion which can be described as `invariance of truncated series sum after correction'. Also Madhava's method of arriving at the correction terms gave rise  to a family of series for $\pi$ !. The method used by Madhava is as follows. In the series 
for $\pi$ let the $m^{th}$ term be $(-1)^m\frac{1}{2m+1}$ let $p=2m+1$
and let the correction applied be $\frac{1}{a_p}$. The criterion used 
is that the series sum should not change
due to correction after $m-1$ or $m$ terms. Thus:
\begin{eqnarray*}
1-\frac{1}{3}+\cdots+(-1)^\frac{p-1}{2}\frac{1}{p} + (-1)^\frac{p+1}{2}\frac{1}{a_p} &=& 1-\frac{1}{3}+\cdots+(-1)^\frac{p-3}{2}\frac{1}{p-1} +  (-1)^\frac{p-1}{2}\frac{1}{a_{p-1}} \\
\frac{1}{a_{p-1}}+\frac{1}{a_p}&=& \frac{1}{p} \\
\end{eqnarray*}
The error due to choice of $a_p$ is captured by the term $\frac{1}{a_{p-1}}+\frac{1}{a_p}-\frac{1}{p}$. Suppose we choose the rule $a_p=2p$ the error is: 
\begin{eqnarray*}
	E(p)&=&\frac{1}{2p-4}+\frac{1}{2p}-\frac{1}{p} \\
	    &=&\frac{1}{(p-1)^2-1} \\
\end{eqnarray*}  
By adding and subtracting $\frac{1}{a_p}$ the original series can be rewritten to  
get a new corrected series for $\pi$ which is \cite[Section~ 6.9-10]{Jyesthadeva2008}:
\begin{eqnarray*}
	\frac{\pi}{4}=\frac{1}{2} + \frac{1}{(2^2-1)} - \frac{1}{(4^2-1)} + \cdots \\
\end{eqnarray*}
It is shown that $a_p=2(p+1)$ is a better choice and that the correction terms can be refined using continued fraction expansion. This gives rise to a sequence of correction terms one better than the other i.e
$\frac{1}{4n}$, $\frac{n}{4n^2+1}$, $\frac{n^2+1}{n(4n^2+5)}$..

Madhava's method gives a correction term, corresponding accelarated series and a method to better the correction as a package.The criterion used by Madhava to find the correction terms which we refer to as `invariance of the truncated series sum after correction' is quite unique. Centuries later many other correction strategies have been proposed for accelaration of series by the likes of Euler, Stirling, Richardson, Aitken, Levin \cite{Corless2023} but Madhava's method is different from all these and deserves much more attention and closer analysis.  
 At the very least this shows Madhava's ability to estimate the reminder of the series with arbitrary accuracy. Thus given
the modern definitions of convergence involving $\epsilon$ and
$N$ the chance that Madhava would have come up with the $N$'s for the $\epsilon$ challenges seems pretty high. And mainly this is what has gone into proving the convergence of the sum of arc-bit approximations in Lemma 2.  

\subsection{Series for $\mathbf{arctan(x),sine(x),cosine(x)}$}
\label{subsec:other_series}
As explained in Section \ref{sec:madhava_other} Madhava's method is general enough to generate other trigonometric series such as the sine series and the cosine series with appropriate changes to 
bit addition, iterative refinement  and the approximation
of power series sums. Leibniz's first step gives a way to 
find series for a function whose differential is a rational function like $\frac{1}{1+x^2}$ which can be expanded as a geometric series. This gives series for functions like ln(1+x), arctanh(x) \cite[Section~3.1,3.11]{stewart2008calculus}. For sine and cosine the popular modern way is to use the Taylor's series
\cite[Section~11.10]{stewart2008calculus}.

\label{sec:Discuss}

\section{Conclusion}
\label{sec:conclusion}
The $\pi$-series is named after Madhava, Gregory and Leibniz based on the chronology of its discovery. The popular modern proof for the series based on Leibniz's ideas uses rigorous definitions and concepts from Calculus many of which were developed
after Leibniz. We extend the same benifit to Madhava's \emph{upapatti} i.e proof
given in \emph{yuktibhasha}:
\begin{enumerate}
\item By using the trigonometric inequality $tan(x)\ge x$ in \emph{capikaranam,
 Making the Arc} Lemma, we show that Madhava's addition of arc-bit approximations($\hat{A}$) gives the desired arc-bit sum $\overset{\frown}{C}$. 
\item By showing in \emph{phala parampara, Series of Results} Lemma that the reminder term of geometric-series can be made as small
as needed for the series to converge to get the $\pi$-series.  
\item By keeping track of the terms to be neglected using big $O$ notation for the approximation of sum of powers, in \emph{samaghata samkalita, Equal Power Sum} Lemma. 
\end{enumerate}
The main steps in the method used by Madhava
i.e \begin{enumerate}
	\item Setting up the `bit sum',
	\item Iterative refinement and
	\item Approximation of power sums
\end{enumerate}
give the arctan, sine and cosine series as well. We compare the details of the steps in case of the $\pi$-series with those of the sine series. 
We present the popular modern day proof for the $\pi$-series which is based on Leibniz's ideas highlighting the first step which is based on Leibniz's idea of `Transmutation' which eventually took the form of the Fundamental Theorem of Calculus (FTC). We compare the two proofs and find that:
\begin{enumerate}
	\item The two proofs are fundamentally different as Madhava's is based on summing the arc-bits where as Leibniz's is based on finding the area under the arc.
	\item The popular proof attributed to Leibniz relies on the machinery of calculus i.e Continuity, Differentiation, Integration upto FTC developed over a period of two centuries starting with Gregory, Leibniz and Newton and ending with Cauchy, Weierstrass and Reymond.
	\item  Madhava's derivation stands out for its simplicity since it can be derived using results related to convergence of series applied specially to the alternating series.
	\item Madhava's series comes with a unique method to add correction terms to accelarate convergence.
\end{enumerate}
 Madhava's proof deserves to be appreciated, taught and popularized as the original proof for the $\pi$-series.

\bibliographystyle{unsrt}
\bibliography{references} 

\end{document}